\documentclass[11pt,a4paper]{article}
\usepackage[margin=2.6cm]{geometry}
\newenvironment{keywords}{\par\noindent\textbf{KEYWORDS}\par\noindent}{\par\medskip}
\newenvironment{amscode}{\par\noindent\textbf{AMS CLASSIFICATION}\par\noindent}{\par}

\usepackage{amsmath,amssymb,amsthm,mathtools}
\usepackage{mathrsfs}
\usepackage{enumitem}

\usepackage[numbers,sort&compress]{natbib}
\bibpunct[, ]{[}{]}{,}{n}{,}{,}

\theoremstyle{plain}
\newtheorem{theorem}{Theorem}[section]
\newtheorem{proposition}[theorem]{Proposition}
\newtheorem{lemma}[theorem]{Lemma}

\theoremstyle{definition}
\newtheorem{definition}[theorem]{Definition}

\theoremstyle{remark}
\newtheorem{remark}[theorem]{Remark}

\newcommand{\R}{\mathbb R}

\newcommand{\ddc}{dd^c}
\newcommand{\SH}{\mathrm{SH}}
\newcommand{\PSH}{\mathrm{PSH}}
\newcommand{\mCap}{\mathrm{Cap}}

\newcommand{\osc}{\mathrm{osc}}

\newcommand{\cP}{\mathscr P}

\newcommand{\one}{\mathbf 1}

\newcommand{\Hm}{H_m}

\newcommand{\norm}[1]{\left\lVert #1\right\rVert}
\newcommand{\abs}[1]{\left|#1\right|}
\newcommand{\set}[1]{\left\{#1\right\}}

\begin{document}

\title{Hessian integration functionals on compact K\"ahler manifolds}

\author{Nguyen Quang Dieu and Pham Le Khanh\thanks{corresponding author}}
\date{}

\maketitle

\begin{abstract}
We show that if $\varphi$ is a H\"older continuous
$\omega$-$m$-subharmonic function on a compact K\"ahler manifold, then its
Hessian measure acts H\"older continuously, for the $L^1$ distance, on all
normalized $\omega$-$m$-subharmonic functions.  We deduce energy bounds and
stability for the associated Hessian equations, and show that a H\"older
classical super-potential does not imply finite-energy solvability at
lower Hessian orders.
\end{abstract}


\begin{amscode}
32U05; 32U15; 32W20; 35J60; 53C55
\end{amscode}

\section{Introduction}

Let $\Omega$ be a domain in $\mathbb C^n$ and let
$\beta:=\ddc |z|^2$ be the standard K\"ahler form.  Two of the basic
nonlinear operators in pluripotential theory are the complex Monge--Amp\`ere
operator
\[
  (\ddc u)^n
\]
and, more generally, the complex $m$-Hessian operator
\[
  (\ddc u)^m\wedge\beta^{n-m},\qquad 1\leq m\leq n.
\]
For $m=n$ one recovers the Monge--Amp\`ere operator, whereas for $m=1$ the
operator is, up to a normalizing constant, the Laplacian.  The intermediate values $1<m<n$ form the
complex Hessian theory.  In this case the natural functions are
$m$-subharmonic functions, and the positivity of the operator is governed by
the G\r{a}rding cone.  The local weak theory was developed in particular by
B{\l}ocki \cite{Blocki}, following the smooth elliptic theory of
Caffarelli--Nirenberg--Spruck \cite{CNS}; a variational approach in
$\mathbb C^n$ was developed by Lu \cite{LuVariational}.  The Bedford--Taylor construction
and capacity theory at the Monge--Amp\`ere endpoint remain an important model
for this development; see \cite{BedfordTaylorCapacity}.

The same equations arise naturally on a K\"ahler manifold.  Let
$(X,\omega)$ be a compact connected K\"ahler manifold of complex dimension
$n$.  If $U\subset X$ is a coordinate ball, then one may write
$\omega=\ddc\rho$ on $U$, where $\rho$ is a smooth strictly
plurisubharmonic function.  An upper semicontinuous function $u$ is
$\omega$-$m$-subharmonic precisely when, locally, $\rho+u$ is
$m$-subharmonic with respect to the background form $\omega|_U$.  On such a coordinate ball one has
\[
  \omega+\ddc u=\ddc(\rho+u)
\]
and hence
\[
  (\omega+\ddc u)^m\wedge\omega^{n-m}
  =
  \bigl(\ddc(\rho+u)\bigr)^m\wedge(\ddc\rho)^{n-m}.
\]
The operator is therefore locally a Hessian operator with a smooth
K\"ahler background.  The definition and approximation results used below
are recalled in Section~\ref{sec:hessian-preliminaries}; see \cite{LuNguyen}.

There are, however, several global features which do not occur in the same
form on a bounded domain.  A compact manifold has no boundary on which one
can prescribe Dirichlet data.  The background form $\omega$ need not admit a
global potential, and the potentials are determined only up to an additive
constant.  On the other hand, closedness of $\omega$ gives an important
global normalization.  If $u$ is bounded and $\omega$-$m$-subharmonic, then
Stokes' theorem gives
\[
  \int_X(\omega+\ddc u)^m\wedge\omega^{n-m}
  =\int_X\omega^n=:V_\omega.
\]
Consequently, on a compact K\"ahler manifold the analogue of a Dirichlet
problem is often a prescribed measure problem
\[
  (\omega+\ddc u)^m\wedge\omega^{n-m}=\mu,
  \qquad \mu(X)=V_\omega,
\]
together with a normalization such as $\sup_Xu=0$.  When singular potentials
are allowed, the finite-energy classes provide the natural replacement for
bounded solutions.  

The Monge--Amp\`ere case gives the clearest example of this global formulation.  The
pluripotential approach to degenerate complex Monge--Amp\`ere equations on
compact K\"ahler manifolds has produced a strong existence and regularity
theory; see, among others, \cite{KolodziejActa,KolodziejHolder}.  For complex
Hessian equations, the corresponding compact theory was developed through
a priori estimates, mixed Hessian inequalities, approximation, capacity and
variational methods; see
\cite{HouMaWu,DinewKolodziej,LuDegenerate,DinewLu,LuNguyen}.  

In this paper we consider a question of this type which originates at the
Monge--Amp\`ere endpoint.  For a positive measure $\mu$ on $X$, one may study
$\mu$ through its action on normalized potentials.  In the Monge--Amp\`ere
case, Dinh and Nguyen \cite{DinhNguyenCMAPH} proved that a positive measure
$\mu$ of mass $V_\omega$ is of the form
\[
  \mu=(\omega+\ddc\varphi)^n
\]
for a H\"older continuous $\omega$-plurisubharmonic potential $\varphi$ if
and only if the corresponding integration functional is H\"older continuous
on normalized test functions that are $\omega$-plurisubharmonic.  

It is natural to ask whether the same mechanism survives for the complex
Hessian operator.  Fix $1\leq m\leq n$ and write
\[
  \Hm(u):=(\omega+\ddc u)^m\wedge\omega^{n-m}.
\]
We use the normalized class
\[
  \cP_m(X,\omega):=
  \{v\in\SH_m(X,\omega):\sup_Xv=0\}
\]
and associate with a finite positive measure $\mu$ the functional
\[
  \mathcal U_\mu^{(m)}:\cP_m(X,\omega)\longrightarrow
  \mathbb R\cup\{-\infty\},
  \qquad
  \mathcal U_\mu^{(m)}(v):=\int_Xv\,d\mu.
\]
Since $\cP_m(X,\omega)$ is relatively compact and uniformly bounded in
$L^1(X,\omega^n)$, the $L^1(\omega^n)$ distance is a natural weak distance on
this class.  We therefore ask when $\mathcal U_\mu^{(m)}$ is H\"older
continuous with respect to this distance.  For $m=n$ this is the potential
form of the classical super-potential considered by Dinh--Nguyen.  For
$m<n$ it tests the measure on a strictly larger family of potentials.

When $m<n$, a test function in $\cP_m(X,\omega)$ need not be
$\omega$-plurisubharmonic, so its current $\omega+\ddc u$ need not be
positive.  Our proof follows \cite[Proposition~3.6]{DinhNguyenCMAPH},
using polarized G\r{a}rding inequalities for the mixed products and the
approximation theory of Lu--Nguyen.  The cohomological mass identities
make the estimates independent of the sup-norms of the bounded tests;
truncation then gives the estimate on the full normalized class.

There are close local antecedents.  Ko{\l}odziej and Nguyen proved that the
Hessian measure of a H\"older continuous $m$-subharmonic function is
quantitatively dominated by the relative Hessian capacity
\cite{KolodziejNguyenCapacity}, and Benali and Zeriahi used a closely
related regularization-and-induction argument for the H\"older continuous
subsolution problem on bounded strongly $m$-pseudoconvex domains
\cite{BenaliZeriahi2020}; see also the erratum \cite{BenaliZeriahiErratum}.
For $m=n$, the qualitative conclusion is contained in
\cite[Lemma~3.5 and Proposition~3.6]{DinhNguyenCMAPH}.  Here we extend it
to the full Hessian test space for $m<n$.

Our first result is the forward estimate for Hessian measures generated by a
H\"older continuous $\omega$-$m$-subharmonic potential.
\begin{theorem}\label{thm:main-necessary}
Let $(X,\omega)$ be a compact K\"ahler manifold of dimension $n$ and let
$1\leq m\leq n$.  Suppose that
\[
 \varphi\in\SH_m(X,\omega)\cap C^{0,\alpha}(X),
 \qquad 0<\alpha\leq1,
\]
and put
\[
 \mu_\varphi:=\Hm(\varphi)
 = (\omega+\ddc\varphi)^m\wedge\omega^{n-m}.
\]
Then every $u\in\cP_m(X,\omega)$ belongs to $L^1(\mu_\varphi)$.  Moreover,
there is a constant $C_m>0$ such that
\[
 \int_X|u-v|\,d\mu_\varphi
 \leq C_m\|u-v\|_{L^1(\omega^n)}^{(\alpha/2)^m},
 \qquad u,v\in\cP_m(X,\omega).
\]
If $\sup_X\varphi=0$ and
$\|\varphi\|_{C^{0,\alpha}(X)}\leq M_\varphi$, then $C_m$ may be chosen to
depend only on $(X,\omega)$, $m$, $\alpha$, and $M_\varphi$; in particular,
it is independent of $u$ and $v$.
\end{theorem}

In particular, $\mathcal U_{\mu_\varphi}^{(m)}$ is
$(\alpha/2)^m$-H\"older continuous for the $L^1(\omega^n)$ distance.
We next record the resulting energy bound and stability statement.
The class $\mathcal E^1(X,\omega,m)$ and its energy $\mathcal E_m$ are
recalled in Section~\ref{sec:finite-energy-solvability}.

\begin{theorem}
\label{thm:weighted-solvability}
Under the hypotheses of Theorem~\ref{thm:main-necessary}, let
\[
 \mu_\varphi:=H_m(\varphi)
\]
and let $0\leq f\in L^\infty(\mu_\varphi)$ satisfy the compatibility condition
\[
 \int_X f\,d\mu_\varphi=V_\omega.
\]
Then the following assertions hold.
\begin{itemize}[leftmargin=*,itemsep=.45\baselineskip]
\item There exists a unique normalized solution
\[
 U\in\mathcal E^1(X,\omega,m),\qquad
 \sup_XU=0,\qquad H_m(U)=f\mu_\varphi.
\]

\item There is a constant $C_{\rm en}>0$, independent of $f$, such that
\[
 -\mathcal E_m(U)
 \leq C_{\rm en}\|f\|_{L^\infty(\mu_\varphi)}.
\]
If $\sup_X\varphi=0$ and
$\|\varphi\|_{C^{0,\alpha}(X)}\leq M_\varphi$, then $C_{\rm en}$ depends only
on $(X,\omega)$, $m$, $\alpha$, and $M_\varphi$.

\item Let $0\leq f_j\in L^\infty(\mu_\varphi)$ satisfy
\[
 \int_Xf_j\,d\mu_\varphi=V_\omega,
 \qquad
 \sup_j\|f_j\|_{L^\infty(\mu_\varphi)}<+\infty,
 \qquad f_j\longrightarrow f\quad\text{in }L^1(\mu_\varphi),
\]
and let $U_j$ be the normalized solution of
$H_m(U_j)=f_j\mu_\varphi$.  Then
\[
 U_j\longrightarrow U\quad\text{in }L^1(\omega^n),
 \qquad
 \mathcal E_m(U_j)\longrightarrow\mathcal E_m(U).
\]
\end{itemize}
\end{theorem}

Existence follows from the range characterization
\cite[Theorem~7.2]{LuNguyen}: since $\varphi$ is bounded, it belongs to
$\mathcal E^1(X,\omega,m)$, so the ``only if'' direction of that theorem
gives $\mathcal E^1(X,\omega,m)\subset L^1(\mu_\varphi)$, and
multiplication by $f\in L^\infty(\mu_\varphi)$ preserves this
integrability.  The functional estimate gives the energy bound and stability.
Remark~\ref{rem:weaker-hypotheses} shows that these conclusions only require
boundedness and continuity of $\varphi$, respectively.  For related
bounded-subsolution results on manifolds with boundary, see
\cite[Theorems~8.7 and~9.1]{KolodziejNguyen} and
\cite{GuNguyen,CollinsPicard}.

Theorem~\ref{thm:weighted-solvability} is a same-order statement.  The next
result shows that a H\"older classical super-potential, which contains only
top-degree information, need not imply finite-energy solvability for a
strictly lower Hessian order.  The construction takes place on
$(\mathbb P^1)^n$ and works for every $1\leq m<n$.
\begin{proposition}
\label{prop:projective-obstruction}
Let $n\geq2$, let $1\leq m<n$, and fix
\[
 \frac{n-m}{m+1}\leq\lambda<\frac{n-m}{m}.
\]
On $X=(\mathbb P^1)^n$, there exist a K\"ahler form $\omega$ cohomologous to
the product Fubini--Study form, a point $x_0\in X$, and a positive measure
$\mu$ of mass $V_\omega$ such that:
\begin{itemize}[leftmargin=*,itemsep=.45\baselineskip]
\item there are $\alpha>0$ and
$\Phi\in\PSH(X,\omega)\cap C^{0,\alpha}(X)$ satisfying
\[
 \mu=(\omega+\ddc\Phi)^n;
\]
in particular, the classical super-potential of $\mu$ is H\"older continuous;

\item the equation $H_m(V)=\mu$ has a unique normalized full-mass solution
$V\in\mathcal E(X,\omega,m)$, and in flat coordinates centered at $x_0$,
\[
 V(z)=-\delta|z|^{-2\lambda}
 \quad\text{near }x_0
\]
for some $\delta>0$;

\item $V\notin\mathcal E^1(X,\omega,m)$.  Consequently, the equation has no
finite-energy solution and no bounded solution.
\end{itemize}
\end{proposition}

The radial model behind Proposition~\ref{prop:projective-obstruction} is of
the type commonly used to test integrability thresholds of $m$-subharmonic
functions; compare \cite{Blocki,DinewKolodziej}.  Likewise, the H\"older
continuity of $\Phi$ follows from Ko{\l}odziej's theorem, since the density
of $\mu$ lies in $L^p$ for some $p>1$.  The point of the proposition is the
cross-order comparison: the same measure is regular from the
Monge--Amp\`ere point of view and singular from the point of view of a lower
Hessian order.

The paper is organized as follows.  Section~\ref{sec:hessian-preliminaries}
recalls the positivity, compactness, and bounded Hessian products used
below.
Section~\ref{sec:hessian-functionals} briefly relates the classical
Monge--Amp\`ere and Hessian test spaces, and develops the elementary
properties of Hessian integration functionals.
Section~\ref{sec:main-proof} contains the proof of
Theorem~\ref{thm:main-necessary}, centered on
Proposition~\ref{prop:inductive}, and then applies the resulting functional
estimate to prove Theorem~\ref{thm:weighted-solvability}.
Section~\ref{sec:consequences} proves
Proposition~\ref{prop:projective-obstruction} and records open questions.

\section{Complex Hessian preliminaries}
\label{sec:hessian-preliminaries}

Throughout this section $(X,\omega)$ is a compact connected K\"ahler
manifold of complex dimension $n$.  We use
\(\ddc=i\partial\bar\partial\) and write
\[
 V_\omega:=\int_X\omega^n.
\]
A different standard normalization of $\ddc$ changes only the numerical
constants.  We first recall the finite-dimensional positivity underlying
the Hessian operator and then the global potential theory needed in the
proofs.

\subsection{The G\r{a}rding cones}

For $\lambda=(\lambda_1,\ldots,\lambda_n)\in\R^n$, let
\[
 \sigma_k(\lambda)
 :=\sum_{1\leq i_1<\cdots<i_k\leq n}
   \lambda_{i_1}\cdots\lambda_{i_k},
 \qquad 1\leq k\leq n.
\]
The open G\r{a}rding cone is
\[
 \Gamma_m:=\set{\lambda\in\R^n:
                 \sigma_1(\lambda)>0,\ldots,\sigma_m(\lambda)>0}.
\]
Its closure is denoted by $\overline{\Gamma}_m$.  If $\gamma$ is a real
$(1,1)$-form and $\omega$ is positive definite, we say that $\gamma$ is
$m$-positive with respect to $\omega$ if its $\omega$-eigenvalues belong to
$\overline{\Gamma}_m$.  Equivalently,
\[
 \gamma^k\wedge\omega^{n-k}\geq0,\qquad 1\leq k\leq m.
\]

The cone $\Gamma_m$ is convex.  Moreover, if
$\gamma_1,\ldots,\gamma_k$ are $m$-positive and $k\leq m$, then
\begin{equation}\label{eq:polarized-garding}
 \gamma_1\wedge\cdots\wedge\gamma_k\wedge\omega^{n-k}
 \geq0.
\end{equation}
This mixed positivity is a consequence of G\r{a}rding's polarization
inequalities \cite{Garding}.  It is the basic positivity property used
throughout the paper.

\subsection{Compact K\"ahler Hessian potentials}

On a coordinate ball $U$ choose a smooth local potential $\rho$ with
$\omega=\ddc\rho$.

\begin{definition}\label{def:compact-msh}
An upper semicontinuous function
$u:X\to\R\cup\{-\infty\}$, not identically $-\infty$, is
$\omega$-$m$-subharmonic if $\rho+u$ is $m$-subharmonic on every such
coordinate ball, where the ball is regarded as the K\"ahler manifold
$(U,\omega|_U)$ and $m$-positivity is computed relative to $\omega|_U$.
The resulting class is denoted by $\SH_m(X,\omega)$, and we write
\[
 \omega_u:=\omega+\ddc u.
\]
\end{definition}

The definition is independent of the local potential because two local
potentials differ by a pluriharmonic function; this is
\cite[Definition~2.5]{LuNguyen}.  The following equivalent current criterion
is the compact form of \cite[Lemma~2.2 and Definition~2.5]{LuNguyen}:
\[
 \omega_u\wedge\gamma_1\wedge\cdots\wedge\gamma_{m-1}
       \wedge\omega^{n-m}\geq0
\]
in the sense of currents for all smooth $m$-positive forms
$\gamma_1,\ldots,\gamma_{m-1}$.  In particular, the first trace inequality
\[
 \omega_u\wedge\omega^{n-1}\geq0
\]
shows that every $\omega$-$m$-subharmonic function is quasi-subharmonic with
a lower Laplacian bound depending only on $\omega$.

For orientation, the local pluripotential theory of $m$-subharmonic functions
and bounded Hessian products originates with B{\l}ocki
\cite{Blocki}.  In the compact setting, an approximation class carrying the
bounded Hessian operator was introduced in \cite{LuDegenerate}.
Lu--Nguyen's approximation theorem \cite[Corollary~3.3]{LuNguyen} proves
that this auxiliary class equals the intrinsic class $\SH_m(X,\omega)$ and
identifies its elements with decreasing limits of smooth
$m$-subharmonic potentials.  Their Section~2.3 recalls the bounded Hessian
products and quasi-continuity on the approximation class, while
\cite[Theorem~4.6]{LuNguyen} develops the corresponding capacity theory.
We also use the mixed inequalities of Dinew--Lu \cite{DinewLu}.

When $m=n$, Definition~\ref{def:compact-msh} gives
$\PSH(X,\omega)$.  The inclusions
\[
 \PSH(X,\omega)=\SH_n(X,\omega)
 \subset\SH_{n-1}(X,\omega)\subset\cdots\subset\SH_1(X,\omega)
\]
follow immediately from the nesting of the G\r{a}rding cones.  The class
$\SH_m(X,\omega)$ is convex and stable under finite maxima.

If $u_1,\ldots,u_k\in\SH_m(X,\omega)\cap L^\infty(X)$ and $k\leq m$, the
Bedford--Taylor construction defines the positive measure
\[
 \omega_{u_1}\wedge\cdots\wedge\omega_{u_k}\wedge\omega^{n-k}.
\]
The construction is local and is compatible with decreasing bounded
approximations.  It is symmetric and multilinear in the closed-current
factors $\omega_{u_j}$.  These facts, including monotone convergence, are
proved in the compact K\"ahler setting in
\cite[Section~2.3]{LuNguyen}.  We use in particular the global decreasing
smooth approximation from \cite[Corollary~3.3]{LuNguyen}; see also
\cite{BedfordTaylorCapacity,Blocki,DinewLu}.
Closedness and Stokes' theorem give the mass identity
\[
 \int_X
 \omega_{u_1}\wedge\cdots\wedge\omega_{u_k}\wedge\omega^{n-k}
 =\int_X\omega^n.
\]
For a single bounded potential $u$ we use the notation
\[
 \Hm(u):=\omega_u^m\wedge\omega^{n-m}.
\]

\subsection{Normalized compact families}

Set
\[
 \cP_m(X,\omega)
 :=\set{u\in\SH_m(X,\omega):\sup_Xu=0}.
\]

\begin{lemma}\label{lem:L1-compactness}
The family $\cP_m(X,\omega)$ is relatively compact in $L^1(X,\omega^n)$.
Moreover,
\[
 \sup_{u\in\cP_m(X,\omega)}\norm{u}_{L^1(\omega^n)}<+\infty.
\]
\end{lemma}

\begin{proof}
The trace inequality
\[
 0\leq\omega_u\wedge\omega^{n-1}
 =\left(1+\frac{\Delta_\omega u}{n}\right)\omega^n
\]
gives $\Delta_\omega u\geq-n$.  Thus $\cP_m(X,\omega)$ is contained in the
normalized compact family of quasi-subharmonic functions with this fixed
Laplacian lower bound.  The standard Hartogs compactness theorem gives both
assertions; see \cite[Chapter~III]{DemaillyBook}.  In the Hessian setting, the
same uniform $L^1$ estimate is used in the proof of
\cite[Theorem~4.10]{LuNguyen}.  Equivalently, the assertions follow at once
from Green representation: the measures
$(\Delta_\omega u+n)\omega^n$ have fixed mass, while $\sup_Xu=0$ controls the
additive constants.
\end{proof}

\begin{lemma}\label{lem:max-normalized}
If $u,v\in\cP_m(X,\omega)$, then $\max\{u,v\}\in\cP_m(X,\omega)$ and
\[
 \norm{\max\{u,v\}-u}_{L^1(\omega^n)}
 \leq\norm{u-v}_{L^1(\omega^n)}.
\]
\end{lemma}

\begin{proof}
The stability under finite maxima gives the first assertion.  Since
$u,v\leq0$ and both have supremum zero, the maximum also has supremum zero.
The second assertion is pointwise.
\end{proof}

\subsection{Hessian capacity and quasi-continuity}

For a Borel set $E\subset X$, define
\[
 \mCap_{\omega,m}(E)
 :=
 \sup_{\substack{v\in\SH_m(X,\omega)\\-1\leq v\leq0}}
 \int_E\Hm(v)
\]
and its outer regularization by
\[
 \mCap_{\omega,m}^*(E)
 :=
 \inf_{\substack{E\subset G\\G\ {\rm open}}}
 \mCap_{\omega,m}(G).
\]
This is the capacity introduced in \cite[Section~2.3]{LuNguyen}.  Its outer
regularization is a Choquet capacity, and the inner and outer capacities
coincide on Borel sets by \cite[Theorem~4.6]{LuNguyen}.  The decreasing
approximation theorem \cite[Corollary~3.3]{LuNguyen}, together with the
definition of the compact Hessian class, gives quasi-continuity of bounded
$\omega$-$m$-subharmonic functions.  This capacity is not used in the proofs below; we recall it because it
enters question (Q3) in Section~\ref{sec:open-questions}.

\section{Hessian integration functionals}
\label{sec:hessian-functionals}

\subsection{Relation with classical super-potentials}
\label{sec:classical-superpotentials}

We only need the potential formulation of the classical theory.  Put
\[
 \mathcal Q_\omega
 :=\left\{v\in\PSH(X,\omega):\int_Xv\,\omega^n=0\right\}.
\]
For a finite positive measure $\mu$, the classical super-potential in the
measure case is represented by
\[
 \mathcal U_\mu:\mathcal Q_\omega\longrightarrow
 \mathbb R\cup\{-\infty\},
 \qquad \mathcal U_\mu(v):=\int_Xv\,d\mu.
\]
This is the top-degree specialization of the theory introduced by Dinh and
Sibony \cite{DinhSibonyActa}.  The comparison of the current metrics with
potential $L^1$ distances in
\cite[Propositions~3.1--3.2]{DinhNguyenCMAPH}, followed by the maximum
reduction in \cite[Lemma~3.5]{DinhNguyenCMAPH}, shows that H\"older
continuity of $\mathcal U_\mu$ is equivalent, up to changing the exponent
and the constant, to an estimate
\[
 \int_X|u-v|\,d\mu
 \leq C\|u-v\|_{L^1(\omega^n)}^\beta,
 \qquad u,v\in\mathcal Q_\omega,
\]
for some $C>0$ and $\beta\in(0,1]$.  Here one uses the uniform
$L^1(\omega^n)$ diameter of the normalized family to absorb the large-distance
part of the estimate in that lemma.

When $m=n$, changing from mean normalization to supremum normalization
identifies this potential picture with the test class $\cP_n(X,\omega)$.  If
$m<n$, however, a function in $\cP_m(X,\omega)$ need not be
$\omega$-plurisubharmonic; equivalently, $\omega+\ddc v$ need not be a
positive current.  Thus the classical super-potential does not act on all
Hessian test potentials.  This is the only distinction from the classical
theory needed here and motivates the following definition.

\subsection{Definition and elementary properties}

\begin{definition}\label{def:m-functional}
Let $\mu$ be a finite positive Radon measure on a compact K\"ahler manifold
$(X,\omega)$.
We say that $\mu$ has a H\"older continuous $m$-integration functional with
exponent $\beta\in(0,1]$ if every $u\in\cP_m(X,\omega)$ is $\mu$-integrable
and there is a constant $C>0$ such that
\begin{equation}\label{eq:signed-holder}
 \abs{\int_X(u-v)\,d\mu}
 \leq C\norm{u-v}_{L^1(\omega^n)}^\beta
\end{equation}
for all $u,v\in\cP_m(X,\omega)$.
\end{definition}

For bounded tests, the signed estimate \eqref{eq:signed-holder} already
controls absolute differences, at the cost of a factor $2$.
\begin{lemma}\label{lem:absolute-criterion}
Assume that \eqref{eq:signed-holder} holds for bounded elements of
$\cP_m(X,\omega)$.  Then
\begin{equation}\label{eq:absolute-holder}
 \int_X\abs{u-v}\,d\mu
 \leq 2C\norm{u-v}_{L^1(\omega^n)}^\beta
\end{equation}
for all bounded $u,v\in\cP_m(X,\omega)$.
\end{lemma}

\begin{proof}
This is exactly the maximum device used in the proof of
\cite[Lemma~3.5]{DinhNguyenCMAPH}, now applied to the normalized Hessian
class, which is stable under maxima.
Put $w:=\max\{u,v\}$.
By Lemma~\ref{lem:max-normalized},
$w\in\cP_m(X,\omega)$.  Since
\[
 \abs{u-v}=(w-u)+(w-v),
\]
we obtain
\begin{align*}
 \int_X\abs{u-v}\,d\mu
 &=\int_X(w-u)\,d\mu+\int_X(w-v)\,d\mu\\
 &\leq C\norm{w-u}_{L^1(\omega^n)}^\beta
       +C\norm{w-v}_{L^1(\omega^n)}^\beta\\
 &\leq2C\norm{u-v}_{L^1(\omega^n)}^\beta,
\end{align*}
where the last step uses Lemma~\ref{lem:max-normalized} and the
monotonicity of $t\mapsto t^\beta$.  This proves \eqref{eq:absolute-holder}.
\end{proof}

Once the constant in \eqref{eq:absolute-holder} does not depend on the
sup-norms of the bounded tests, truncation extends the estimate to the whole
normalized class.
\begin{lemma}\label{lem:truncation-extension}
Suppose that \eqref{eq:absolute-holder} holds for all bounded
$u,v\in\cP_m(X,\omega)$, with a constant independent of their sup-norms.
Then every $u\in\cP_m(X,\omega)$ is $\mu$-integrable and the same estimate
holds for arbitrary $u,v\in\cP_m(X,\omega)$.
\end{lemma}

\begin{proof}
For $j\geq1$, set $u_j:=\max\{u,-j\}$.  Then
$u_j\in\cP_m(X,\omega)$, $u_j\searrow u$, and
$u_j\to u$ in $L^1(\omega^n)$.  The assumed estimate shows that
$(u_j)$ is Cauchy in $L^1(\mu)$.  A subsequence converges $\mu$-almost
everywhere to its $L^1(\mu)$ limit.  Since the full sequence decreases
pointwise to $u$, that limit agrees with $u$ $\mu$-almost everywhere.  Thus
$u\in L^1(\mu)$.
Applying the bounded estimate to $u_j,v_j$ and passing to the limit gives the
last assertion.
\end{proof}

\section{H\"older Hessian potentials: proof of the main theorem}
\label{sec:main-proof}

Throughout this section $(X,\omega)$ is a compact K\"ahler manifold, and
\[
 V_\omega:=\int_X\omega^n.
\]
Fix
\[
 \varphi\in\SH_m(X,\omega)\cap C^{0,\alpha}(X),
 \qquad 0<\alpha\leq1,
\]
and normalize it by $\sup_X\varphi=0$.
Let
\[
 A:=\omega+\ddc\varphi,\qquad
 \mu_k:=A^k\wedge\omega^{n-k},\quad 0\leq k\leq m.
\]
In particular, $\mu_0=\omega^n$ and $\mu_m=\Hm(\varphi)$.
The proof follows the architecture of
\cite[Proposition~3.6]{DinhNguyenCMAPH}: pass through the measures $\mu_k$,
reduce arbitrary pairs to ordered ones by maxima, regularize $\varphi$, and
integrate by parts.  What requires care is not the one-variable optimization,
but the geometry behind it.  Since the tests are only
$\omega$-$m$-subharmonic, neither the positivity of the intermediate current
nor the Bedford--Taylor integration-by-parts formula is automatic from the
Monge--Amp\`ere theory.  Lemmas~\ref{lem:mixed-current}
and~\ref{lem:BT-integration-parts} provide precisely these two Hessian inputs.

Once they are available, the larger Hessian cone becomes an advantage: the
same argument applies to every test in $\cP_m(X,\omega)$ and to every mixed
measure between $\omega^n$ and $\Hm(\varphi)$.  We retain the sharper scale
$\|\varphi_\varepsilon\|_{C^2}=O(\varepsilon^{\alpha-2})$ furnished by
H\"older regularization.  This quantitative point gives the recurrence
$\beta_k=(\alpha/2)\beta_{k-1}$.

\subsection{A regularization lemma}

\begin{lemma}\label{lem:holder-regularization}
Let $\varphi\in C^{0,\alpha}(X)$, $0<\alpha\leq1$.  There are smooth functions
$\varphi_\varepsilon$, $0<\varepsilon<\varepsilon_0$, such that
\begin{equation}\label{eq:regularization-bounds}
 \norm{\varphi-\varphi_\varepsilon}_{L^\infty}
 \leq C_{\rm reg}\varepsilon^\alpha,
 \qquad
 -C_{\rm reg}\varepsilon^{\alpha-2}\omega
 \leq\ddc\varphi_\varepsilon
 \leq C_{\rm reg}\varepsilon^{\alpha-2}\omega.
\end{equation}
The constants $C_{\rm reg}$ and $\varepsilon_0$ depend only on the fixed
atlas and metric, on $\alpha$, and on an upper bound for $\|\varphi\|_{C^{0,\alpha}(X)}$.
\end{lemma}

\begin{proof}
This is the standard local smoothing used in the proof of
\cite[Proposition~3.6]{DinhNguyenCMAPH}, with the dependence on the H\"older
exponent retained.  We give the short argument because the sharper bound
$O(\varepsilon^{\alpha-2})$, rather than the coarser
$O(\varepsilon^{-2})$, determines the exponent in
Proposition~\ref{prop:inductive}.
Choose coordinate balls $U_j'\Subset U_j$ covering $X$, a subordinate smooth
partition of unity, and standard radial convolution kernels in the
coordinates on $U_j$.  Convolving at scale $\varepsilon$ and summing with the
partition of unity gives a global smooth function
$\varphi_\varepsilon$.  The local H\"older estimate and differentiation of
the kernel give
\[
 \|\varphi-\varphi_\varepsilon\|_\infty\leq C\varepsilon^\alpha,
 \qquad
 \|\varphi_\varepsilon\|_{C^2}\leq C\varepsilon^{\alpha-2}.
\]
For completeness, the second estimate uses the cancellation
$\int D^2\rho=0$: in local coordinates,
\[
 D^2(\varphi*\rho_\varepsilon)(x)
 =
 \int D^2\rho_\varepsilon(y)
       \bigl(\varphi(x-y)-\varphi(x)\bigr)\,dy,
\]
and the right-hand side is bounded by
$C\varepsilon^{\alpha-2}[\varphi]_{C^{0,\alpha}}$.
The same argument applies to the cutoff pieces.  All coordinate changes and
cutoff terms are controlled uniformly because the covering is finite.  The
$C^2$ estimate implies the two-sided form estimate in
\eqref{eq:regularization-bounds}.
\end{proof}

\subsection{Mixed currents and integration by parts}

The induction requires two facts about the mixed currents
$A^{k-1}\wedge\omega^{n-k}$.  The first is their positivity together with
the cohomological mass identity.
\begin{lemma}\label{lem:mixed-current}
Let $1\leq k\leq m$ and let
$\varphi,u\in\SH_m(X,\omega)\cap L^\infty(X)$.  Put
\[
 A:=\omega_\varphi,
 \qquad
 T_{k-1}:=A^{k-1}\wedge\omega^{n-k}.
\]
Then $T_{k-1}$ is a positive closed $(n-1,n-1)$-current,
$\omega_u\wedge T_{k-1}$ is a positive measure, and
\begin{equation}\label{eq:mixed-current-mass}
 \int_X\omega_u\wedge T_{k-1}=V_\omega.
\end{equation}
\end{lemma}

\begin{proof}
Assume first that $\varphi$ and $u$ are smooth.  Then $A$ and $\omega_u$
are $m$-positive.  To see that the $(n-1,n-1)$-form
$T_{k-1}=A^{k-1}\wedge\omega^{n-k}$ is positive, let $\gamma$ be an
arbitrary smooth positive $(1,1)$-form.  Its eigenvalues are nonnegative,
so $\gamma$ is $m$-positive, and \eqref{eq:polarized-garding}, applied to
$\gamma$ together with $k-1$ copies of $A$ (in total $k\leq m$ factors),
gives
\[
 \gamma\wedge A^{k-1}\wedge\omega^{n-k}\geq0.
\]
Since this holds for every positive $(1,1)$-form $\gamma$, the form
$T_{k-1}$ is positive; in bidegree $(n-1,n-1)$ weak and strong positivity
coincide.  Applying \eqref{eq:polarized-garding} in the same way to
$\omega_u$ together with $k-1$ copies of $A$ gives
\[
 \omega_u\wedge A^{k-1}\wedge\omega^{n-k}\geq0.
\]
All factors are closed and cohomologous to $\omega$.  Hence Stokes' theorem
gives
\[
 \int_X\omega_u\wedge A^{k-1}\wedge\omega^{n-k}
 =\int_X\omega^n=V_\omega.
\]

We now pass to bounded potentials.  By
\cite[Corollary~3.3]{LuNguyen}, choose smooth
$\varphi_j,u_j\in\SH_m(X,\omega)$ with
$\varphi_j\searrow\varphi$ and $u_j\searrow u$.  Since the limits are
bounded, the two approximating families are uniformly bounded.  Set
\[
 A_j:=\omega_{\varphi_j},\qquad
 T_{k-1,j}:=A_j^{k-1}\wedge\omega^{n-k},\qquad
 \nu_j:=\omega_{u_j}\wedge T_{k-1,j}.
\]
By \cite[Corollary~3.3]{LuNguyen}, bounded $\omega$-$m$-subharmonic
functions belong to Lu's approximation class.  The decreasing sequences
above converge quasi-uniformly by
\cite[Proposition~3.20]{LuDegenerate}.  We may therefore apply the
mixed-product convergence theorem
\cite[Proposition~3.12]{LuDegenerate} to the factors (adjoining zero
potentials where necessary).  Thus
\[
 T_{k-1,j}\longrightarrow
 A^{k-1}\wedge\omega^{n-k}=T_{k-1},
 \qquad
 \nu_j\longrightarrow
 \omega_u\wedge A^{k-1}\wedge\omega^{n-k}
\]
weakly as currents and measures.  Each $T_{k-1,j}$ is positive and closed
and each $\nu_j$ is positive, so the same is true of the limits.  Finally,
$\int_X\nu_j=V_\omega$ for every $j$ by the smooth case, and weak
convergence against the constant function $1$ yields
\eqref{eq:mixed-current-mass}.  This also makes explicit that the mass
identity passes to the simultaneous mixed limit.
\end{proof}

The second fact is an integration-by-parts identity in which the difference
of two bounded test functions is paired with the difference of a bounded
$\omega$-$m$-subharmonic function and a smooth function.
\begin{lemma}\label{lem:BT-integration-parts}
Let $T$ be a mixed current of the form
\[
 T=\omega_{w_1}\wedge\cdots\wedge\omega_{w_r}
   \wedge\omega^{n-1-r},
 \qquad 0\leq r\leq m-1,
\]
where $w_1,\ldots,w_r\in\SH_m(X,\omega)\cap L^\infty(X)$.
Let $u,v,\phi\in\SH_m(X,\omega)\cap L^\infty(X)$ and let
$\psi\in C^\infty(X)$.  Then
\begin{equation}\label{eq:BT-integration-parts}
 \int_X(u-v)\,\ddc(\phi-\psi)\wedge T
 =
 \int_X(\phi-\psi)\,\ddc(u-v)\wedge T.
\end{equation}
Here both sides are understood as the finite sums obtained by expanding the
differences into bounded Bedford--Taylor pairings and smooth pairings, as
made explicit in the proof.
\end{lemma}

\begin{proof}
Insert zero potentials if necessary, so that
\[
 T=\omega_{w_1}\wedge\cdots\wedge\omega_{w_{m-1}}
   \wedge\omega^{n-m},
\]
with bounded $w_j\in\SH_m(X,\omega)$.  By
\cite[Corollary~3.3]{LuNguyen}, these functions and the bounded functions
$u,v,\phi$ belong to the approximation class used in
\cite{LuDegenerate}.  The Bedford--Taylor integration-by-parts theorem
\cite[Theorem~3.13]{LuDegenerate} therefore gives
\[
 \int_X u\,\ddc\phi\wedge T
 =\int_X\phi\,\ddc u\wedge T,
 \qquad
 \int_X v\,\ddc\phi\wedge T
 =\int_X\phi\,\ddc v\wedge T.
\]
All these terms are defined by the bounded mixed Hessian calculus; adjoining
one further Hessian factor gives at most $m$ nontrivial factors.

For the smooth function $\psi$, choose $\tau>0$ so small that
\[
 \omega+\tau\ddc\psi\geq\frac12\omega.
\]
Then $\tau\psi$ is $\omega$-plurisubharmonic and hence
$\omega$-$m$-subharmonic.  Applying the same theorem to
$(u,\tau\psi)$ and $(v,\tau\psi)$, and dividing by $\tau$, gives
\[
 \int_X u\,\ddc\psi\wedge T
 =\int_X\psi\,\ddc u\wedge T,
 \qquad
 \int_X v\,\ddc\psi\wedge T
 =\int_X\psi\,\ddc v\wedge T.
\]
Subtracting these four identities proves
\eqref{eq:BT-integration-parts}.  Thus no bilinear pairing on arbitrary
differences of $\omega$-$m$-subharmonic functions is being postulated:
only this finite combination of standard bounded Bedford--Taylor pairings is
used below.
\end{proof}

\subsection{The fundamental inductive estimate}

The following proposition is the main technical result of the paper.  Its
induction is the Hessian counterpart of
\cite[Proposition~3.6]{DinhNguyenCMAPH}, but the statement is not contained
there: it controls the measures on the larger test class
$\cP_m(X,\omega)$, passes through $m$-positive mixed currents, and remains
uniform under truncation of unbounded tests.  It is stronger than
Theorem~\ref{thm:main-necessary} in that it controls every mixed measure
between $\omega^n$ and $\Hm(\varphi)$.  Its formulation also separates the
Hessian input---positivity and mass of the mixed currents---from the
one-variable optimization that produces the H\"older exponent.
When $m=n$, the proof follows the Dinh--Nguyen architecture, with the explicit
regularization scale above and the absolute-value/truncation bookkeeping made
uniform.  The additional step for $m<n$ is the preservation of
positivity and mass after the plurisubharmonic test class is replaced by
$\cP_m(X,\omega)$; this is exactly the role of
Lemmas~\ref{lem:mixed-current} and~\ref{lem:BT-integration-parts}.
We now state and prove the estimate at every intermediate level.
\begin{proposition}\label{prop:inductive}
For $0\leq k\leq m$ there are constants $C_k>0$ and $\beta_k>0$ such that
\[
 \int_X\abs{u-v}\,d\mu_k
 \leq C_k\norm{u-v}_{L^1(\omega^n)}^{\beta_k}
\]
for all $u,v\in\cP_m(X,\omega)$.  One may take
\[
 \beta_0=1,\qquad
 \beta_k=\frac{\alpha}{2}\beta_{k-1}.
\]
In particular,
\[
 \beta_k=\left(\frac{\alpha}{2}\right)^k.
\]
\end{proposition}

\begin{proof}
For $k=0$,
\[
 \int_X|u-v|\,d\mu_0
 =\|u-v\|_{L^1(\omega^n)},
\]
so the assertion holds with $C_0=1$ and $\beta_0=1$.

Fix $1\leq k\leq m$ and assume the assertion at level $k-1$.  We begin with
a bounded ordered pair $u,v\in\cP_m(X,\omega)$, $u\geq v$, and set
\[
 f:=u-v\geq0,\qquad
 T_{k-1}:=A^{k-1}\wedge\omega^{n-k}.
\]
By Lemma~\ref{lem:mixed-current}, $T_{k-1}$ is a positive closed
$(n-1,n-1)$-current.  Moreover,
$\omega_u\wedge T_{k-1}$ and $\omega_v\wedge T_{k-1}$ are positive measures,
each of mass $V_\omega$.  Since
\[
 \mu_k=\mu_{k-1}+\ddc\varphi\wedge T_{k-1},
\]
the regularization in Lemma~\ref{lem:holder-regularization} gives, for
$0<\varepsilon<\varepsilon_0$,
\begin{align*}
 \int_X f\,d\mu_k
 &=
 \int_X f\,d\mu_{k-1}
 +\int_X f\,\ddc\varphi_\varepsilon\wedge T_{k-1}\\
 &\quad
 +\int_Xf\,\ddc(\varphi-\varphi_\varepsilon)
                    \wedge T_{k-1}.
\end{align*}
Since $f\geq0$, the two-sided form estimate for
$\ddc\varphi_\varepsilon$ and the positivity of $T_{k-1}$ imply
\[
 \abs{\int_X f\,\ddc\varphi_\varepsilon\wedge T_{k-1}}
 \leq C\varepsilon^{\alpha-2}
       \int_Xf\,\omega\wedge T_{k-1}
 =C\varepsilon^{\alpha-2}\int_Xf\,d\mu_{k-1}.
\]

Put $h_\varepsilon:=\varphi-\varphi_\varepsilon$ and recall that
$f=u-v$.  Lemma~\ref{lem:BT-integration-parts}, applied to the same $u$ and $v$
with $\phi=\varphi$ and $\psi=\varphi_\varepsilon$, gives
\begin{align*}
 &\int_Xf\,\ddc(\varphi-\varphi_\varepsilon)\wedge T_{k-1}
 \\
 &\qquad
 =\int_X(\varphi-\varphi_\varepsilon)\,
          \ddc(u-v)\wedge T_{k-1}
 \\
 &\qquad
 =\int_X(\varphi-\varphi_\varepsilon)\,
      (\omega_u-\omega_v)\wedge T_{k-1}.
\end{align*}
The two measures in the last line are positive and, by
\eqref{eq:mixed-current-mass}, have the same mass:
\[
 \int_X\omega_u\wedge T_{k-1}
 =
 \int_X\omega_v\wedge T_{k-1}
 =
 \int_X\omega^n
 =V_\omega.
\]
Since $\|h_\varepsilon\|_\infty\leq C\varepsilon^\alpha$, we obtain
\[
 \abs{\int_Xf\,\ddc(\varphi-\varphi_\varepsilon)
               \wedge T_{k-1}}
 \leq2CV_\omega\varepsilon^\alpha.
\]

Combining the preceding decomposition and estimates gives
\[
 \int_Xf\,d\mu_k
 \leq
 \left(1+C\varepsilon^{\alpha-2}\right)
 \int_Xf\,d\mu_{k-1}
 +C\varepsilon^\alpha.
\]
After decreasing $\varepsilon_0$ if necessary, assume
$0<\varepsilon_0<1$.  Since $\alpha-2<0$, we then have
$1\leq\varepsilon^{\alpha-2}$, and therefore
\(1+C\varepsilon^{\alpha-2}\leq C'\varepsilon^{\alpha-2}\).
The induction hypothesis yields, for
$t:=\norm{u-v}_{L^1(\omega^n)}$,
\begin{equation}\label{eq:induction-before-optimization}
 \int_Xf\,d\mu_k
 \leq C\left(
 \varepsilon^{\alpha-2}t^{\beta_{k-1}}+\varepsilon^\alpha
 \right).
\end{equation}

If $t=0$, then $u=v$ as distributions and hence as their upper
semicontinuous representatives.  Assume $t>0$ and put
\[
 q:=\frac{\beta_{k-1}}{2},
 \qquad
 \varepsilon_*:=\min\left\{\frac{\varepsilon_0}{2},\frac12\right\},
 \qquad
 t_*:=\varepsilon_*^{1/q}.
\]
For $0<t<t_*$, take $\varepsilon=t^q$ in
\eqref{eq:induction-before-optimization}.  The two exponents agree because
\[
 \beta_{k-1}+(\alpha-2)q
 =\alpha q
 =\frac{\alpha}{2}\beta_{k-1}
 =:\beta_k.
\]
It follows that
\[
 \int_Xf\,d\mu_k\leq Ct^{\beta_k},
 \qquad 0<t<t_*.
\]

For $t\geq t_*$, take $\varepsilon=\varepsilon_*$ in
\eqref{eq:induction-before-optimization}.  This gives
\[
 \int_Xf\,d\mu_k
 \leq C\left(1+t^{\beta_{k-1}}\right).
\]
By Lemma~\ref{lem:L1-compactness}, there is $M>0$, independent of $u$ and
$v$, such that $t\leq M$.  On the compact interval $[t_*,M]$ one has
\[
 1+t^{\beta_{k-1}}\leq C_*t^{\beta_k}.
\]
This proves the ordered estimate for every $t\geq0$.

For arbitrary bounded $u,v\in\cP_m(X,\omega)$, put
$w:=\max\{u,v\}$.  Applying the ordered estimate to $(w,u)$ and $(w,v)$ and
using Lemma~\ref{lem:max-normalized} gives
\[
 \int_X|u-v|\,d\mu_k
 =\int_X(w-u)\,d\mu_k+\int_X(w-v)\,d\mu_k
 \leq 2C\norm{u-v}_{L^1(\omega^n)}^{\beta_k}.
\]
Thus the bounded estimate holds at level $k$.  Induction proves it for every
$0\leq k\leq m$.  For each fixed $k$, the constants are independent of the
sup-norms of the bounded tests.  Since $\mu_k$ is a finite positive measure,
Lemma~\ref{lem:truncation-extension} shows that every element of
$\cP_m(X,\omega)$ is $\mu_k$-integrable and extends the same estimate to
arbitrary tests.  The recurrence gives
$\beta_k=(\alpha/2)^k$.
\end{proof}

\begin{proof}[Proof of Theorem~\ref{thm:main-necessary}]
The intermediate measures in Proposition~\ref{prop:inductive} are needed
for the induction but not for the statement of the theorem.  Taking the
endpoint $k=m$ in that proposition gives
\[
 \cP_m(X,\omega)\subset L^1(\mu_m),
 \qquad
 \mu_m=(\omega+\ddc\varphi)^m\wedge\omega^{n-m}=\mu_\varphi,
\]
as well as
\[
 \int_X|u-v|\,d\mu_\varphi
 \leq C_m\|u-v\|_{L^1(\omega^n)}^{\beta_m}.
\]
The functional estimate follows from
\[
 \left|\int_X(u-v)\,d\mu_\varphi\right|
 \leq\int_X|u-v|\,d\mu_\varphi,
\]
and the recurrence in Proposition~\ref{prop:inductive} gives
$\beta_m=(\alpha/2)^m$.

It remains only to record the uniformity asserted in the last sentence.
Under the normalization $\sup_X\varphi=0$ and the bound
$\|\varphi\|_{C^{0,\alpha}}\leq M_\varphi$, the constants
$C_{\rm reg}$ and $\varepsilon_0$ in
Lemma~\ref{lem:holder-regularization} depend only on
$(X,\omega),\alpha$, and $M_\varphi$.  The uniform
$L^1(\omega^n)$-diameter
\[
 M_0:=\sup_{u,v\in\cP_m(X,\omega)}
       \|u-v\|_{L^1(\omega^n)}
 <+\infty
\]
depends only on $(X,\omega)$ and $m$.  At the $k$-th induction step, put
\[
 q_k:=\frac{\beta_{k-1}}2,\qquad
 \varepsilon_*:=\min\{\varepsilon_0/2,1/2\},\qquad
 t_{*,k}:=\varepsilon_*^{1/q_k}.
\]
The proof of Proposition~\ref{prop:inductive} shows that, after the
ordered-pair estimate has been reduced to
\[
 C_k'\bigl(
 \varepsilon^{\alpha-2}t^{\beta_{k-1}}
 +\varepsilon^\alpha\bigr),
\]
one may take an ordered-pair constant bounded by
\[
 C_k'
 \max\left\{
 2,\,
 \varepsilon_*^{\alpha-2}\bigl(1+M_0^{\beta_{k-1}}\bigr)t_{*,k}^{-\beta_k}
 \right\}.
\]
Indeed, for $0<t<t_{*,k}$ the choice $\varepsilon=t^{q_k}$ gives the bound
$2C_k't^{\beta_k}$.  For $t_{*,k}\leq t\leq M_0$ the choice
$\varepsilon=\varepsilon_*$, together with
$\varepsilon_*^{\alpha}\leq\varepsilon_*^{\alpha-2}$, gives the bound
$C_k'\varepsilon_*^{\alpha-2}(1+M_0^{\beta_{k-1}})$, which is at most
$C_k'\varepsilon_*^{\alpha-2}(1+M_0^{\beta_{k-1}})t_{*,k}^{-\beta_k}t^{\beta_k}$.
The maximum reduction introduces only the factor $2$.  Starting from
$C_0=1$ and iterating up to $k=m$ proves that $C_m$ depends only on the data
listed in the theorem and is independent of the test functions.
\end{proof}

\subsection{Finite-energy solvability}
\label{sec:finite-energy-solvability}

We recall the finite-energy facts needed for
Theorem~\ref{thm:weighted-solvability}; see
\cite[Sections~5--7]{LuNguyen}.  For bounded
$w\in\SH_m(X,\omega)$, put
\[
 \mathcal E_m(w)
 :=\frac1{m+1}\sum_{k=0}^m
 \int_X w\,\omega_w^k\wedge\omega^{n-k}.
\]
For general $w\in\SH_m(X,\omega)$, this energy is the infimum of
$\mathcal E_m(v)$ over bounded $\omega$-$m$-subharmonic $v\geq w$.
Write
\[
 \mathcal E^1(X,\omega,m)
 :=\{w\in\mathcal E(X,\omega,m):\mathcal E_m(w)>-\infty\}.
\]
The results of \cite{LuNguyen}, stated there for $V_\omega=1$, apply
for arbitrary $V_\omega$ by replacing $\omega,w$ with
$c\omega,cw$, where $c:=V_\omega^{-1/n}$.  Under this change the
Hessian measure scales by $c^n$, while the energy and linear
functionals scale by $c^{n+1}$.

We shall use three consequences.  First, a positive measure $\nu$ of
mass $V_\omega$ has a solution $H_m(U)=\nu$ in $\mathcal E^1$ if and
only if
\begin{equation}\label{eq:variational-range-criterion}
 \mathcal E^1(X,\omega,m)\subset L^1(\nu);
\end{equation}
see \cite[Theorem~7.2]{LuNguyen}.  Second, for $U\in\mathcal E^1$
with $\sup_XU=0$,
\begin{equation}\label{eq:energy-comparison}
 -\mathcal E_m(U)\leq\int_X(-U)\,H_m(U),
\end{equation}
by \cite[Lemma~6.2]{LuNguyen}.  Such a solution is unique up to a
constant by \cite[Theorem~4.1]{DinewLu}.

Finally, for a positive measure $\nu$ of mass $V_\omega$ satisfying
\eqref{eq:variational-range-criterion}, $H_m(U)=\nu$ is equivalent to
$U$ maximizing
\[
 \mathcal F_\nu(v):=\mathcal E_m(v)-\int_Xv\,d\nu
 \quad\text{on }\mathcal E^1(X,\omega,m).
\]
Indeed, one direction follows from the concavity inequality in
\cite[Lemma~6.3]{LuNguyen}.  Conversely, for a maximizer $W$ and
$q\in C(X,\mathbb R)$, the $\omega$-$m$-subharmonic envelope
$P(W+tq)$ differs from $W$ by at most $|t|\|q\|_\infty$.  Maximality
and positivity of $\nu$ show that
$\mathcal E_m(P(W+tq))-\int_X(W+tq)\,d\nu$ has a maximum at $t=0$.
Its derivative there is
$\int_Xq\,(H_m(W)-\nu)$ by
\cite[Theorem~6.14]{LuNguyen}, so $H_m(W)=\nu$.
Since $\nu(X)=V_\omega$, $\mathcal F_\nu$ is invariant under addition
of constants, and maximization may be restricted to $\sup_Xv=0$.
We also use upper semicontinuity of $\mathcal E_m$ and compactness
of normalized sets with a fixed lower energy bound from
\cite[Lemmas~6.4--6.5]{LuNguyen}.

\begin{proof}[Proof of Theorem~\ref{thm:weighted-solvability}]
Put $\mu:=\mu_\varphi$ and $\beta:=(\alpha/2)^m$.  By
Lemma~\ref{lem:L1-compactness} and
Theorem~\ref{thm:main-necessary},
\[
 M_0:=\sup_{v\in\cP_m(X,\omega)}\|v\|_{L^1(\omega^n)}<\infty,
 \qquad
 \sup_{v\in\cP_m(X,\omega)}\int_X(-v)\,d\mu
 \leq C_mM_0^\beta.
\]
If $w\in\mathcal E^1(X,\omega,m)$, then $w-\sup_Xw\in\cP_m(X,\omega)$,
so $w\in L^1(\mu)$.  Thus $\mathcal E^1\subset L^1(f\mu)$ whenever
$f\in L^\infty(\mu)$.
The range criterion \eqref{eq:variational-range-criterion} gives the
unique normalized solution $U$ of $H_m(U)=f\mu$, and
\eqref{eq:energy-comparison} gives
\[
 -\mathcal E_m(U)
 \leq\int_X(-U)f\,d\mu
 \leq C_mM_0^\beta\|f\|_{L^\infty(\mu)}.
\]
This proves the first two assertions.

For stability, let
$B:=\sup_j\|f_j\|_{L^\infty(\mu)}+\|f\|_{L^\infty(\mu)}$.
For any fixed $v\in\cP_m(X,\omega)$, truncate it by
$v^{(T)}:=\max\{v,-T\}$.  Then
\[
 \left|\int_Xv(f_j-f)\,d\mu\right|
 \leq T\|f_j-f\|_{L^1(\mu)}
       +B\|v-v^{(T)}\|_{L^1(\mu)}.
\]
Taking first $j\to\infty$ and then $T\to\infty$ proves pointwise
convergence of these linear functionals.  Theorem~\ref{thm:main-necessary}
transfers total boundedness of $\cP_m(X,\omega)$ from $L^1(\omega^n)$
to $L^1(\mu)$.  Since the functionals are uniformly $B$-Lipschitz
on $L^1(\mu)$, a finite-net argument gives
\begin{equation}\label{eq:uniform-linear-functional-bound}
 \sup_{v\in\cP_m(X,\omega)}
 \left|\int_Xv(f_j-f)\,d\mu\right|\longrightarrow0.
\end{equation}

Let $\mathcal F_j:=\mathcal F_{f_j\mu}$ and
$\mathcal F:=\mathcal F_{f\mu}$.  The energy bound and
\cite[Lemma~6.5]{LuNguyen} give, after passage to a subsequence,
$U_j\to W\in\mathcal E^1$ in $L^1(\omega^n)$, with $\sup_XW=0$.
By \eqref{eq:uniform-linear-functional-bound} and
Theorem~\ref{thm:main-necessary},
\[
 \int_XU_jf_j\,d\mu\longrightarrow\int_XWf\,d\mu.
\]
Upper semicontinuity of the energy and maximality of $U_j$ imply
\[
 \mathcal F(W)\geq\limsup_j\mathcal F_j(U_j)
 \geq\lim_j\mathcal F_j(v)=\mathcal F(v)
\]
for every normalized $v\in\mathcal E^1$.  Hence $H_m(W)=f\mu$,
so uniqueness gives $W=U$.  Every convergent subsequence has this
limit, and therefore $U_j\to U$ in $L^1(\omega^n)$.
Moreover, testing maximality at $U$ gives
\[
 \mathcal F(U)\leq\liminf_j\mathcal F_j(U_j)
 \leq\limsup_j\mathcal F_j(U_j)\leq\mathcal F(U).
\]
The linear terms converge as above, so
$\mathcal E_m(U_j)\to\mathcal E_m(U)$.
\end{proof}

\begin{remark}\label{rem:weaker-hypotheses}
The regularity assumptions in Theorem~\ref{thm:weighted-solvability}
can be weakened.
\begin{enumerate}[label=\textup{(\roman*)},leftmargin=*]
\item For the existence and energy bound it suffices that
$\varphi\in\SH_m(X,\omega)\cap L^\infty(X)$.  Indeed, for bounded
$v\in\cP_m(X,\omega)$, Lemmas~\ref{lem:mixed-current} and
\ref{lem:BT-integration-parts} give
\[
 \int_X(-v)\,d\mu_k
 =\int_X(-v)\,d\mu_{k-1}
  +\int_X\varphi\,(\omega-\omega_v)\wedge T_{k-1}
 \leq\int_X(-v)\,d\mu_{k-1}+V_\omega\osc_X\varphi.
\]
Iteration and truncation yield
\[
 \sup_{v\in\cP_m(X,\omega)}\int_X(-v)\,H_m(\varphi)
 \leq M_0+mV_\omega\osc_X\varphi.
\]
The range criterion and \eqref{eq:energy-comparison} now give the first
two assertions, with this quantity in place of $C_mM_0^\beta$.

\item For stability it suffices that
$\varphi\in\SH_m(X,\omega)\cap C^0(X)$.  Choose smooth
$\varphi_\varepsilon$ such that
$\delta_\varepsilon:=\|\varphi-\varphi_\varepsilon\|_\infty\to0$
and $-c_\varepsilon\omega\leq\ddc\varphi_\varepsilon
\leq c_\varepsilon\omega$, where $c_\varepsilon<\infty$.
For bounded ordered $u\geq v$, the proof of
Proposition~\ref{prop:inductive} gives
\[
 \int_X(u-v)\,d\mu_k
 \leq(1+c_\varepsilon)\int_X(u-v)\,d\mu_{k-1}
       +2V_\omega\delta_\varepsilon.
\]
Set $\theta_0(t):=t$ and define recursively
\[
 \theta_k(t):=2\inf_{0<\varepsilon<\varepsilon_0}
 \bigl[(1+c_\varepsilon)\theta_{k-1}(t)
                  +2V_\omega\delta_\varepsilon\bigr].
\]
By induction each $\theta_k$ is nondecreasing, being an infimum of
nondecreasing functions of $t$.  Hence taking maxima proves
$\int_X|u-v|\,d\mu_k\leq
\theta_k(\|u-v\|_{L^1(\omega^n)})$ for bounded tests.
For each fixed $\varepsilon$, first let $t\downarrow0$ and then
$\varepsilon\downarrow0$.  Induction gives $\theta_k(0^+)=0$.
The truncation argument in the proof of
Lemma~\ref{lem:truncation-extension} extends the estimate to all
normalized tests.  Thus the $L^1(\mu_\varphi)$-compactness and
uniform-linear-functional argument in the proof above remain valid,
which gives stability.
\end{enumerate}
\end{remark}

\section{A cross-order obstruction at every lower Hessian order}
\label{sec:consequences}

\subsection{Proof of Proposition~\ref{prop:projective-obstruction}}
\label{sec:obstruction}

We work first in the larger range
\[
 0<\lambda<\frac{n-m}{m};
\]
the proposition corresponds to
\(\lambda\geq(n-m)/(m+1)\).  The construction has four steps, separating the
global geometry, the local admissibility, the sharp integrability thresholds,
and the top-degree realization of the measure.

\medskip
\noindent\emph{Step 1: a flat coordinate ball.}
Let \(X=(\mathbb P^1)^n\) and
\[
 \omega_0=\sum_{j=1}^n\operatorname{pr}_j^*\omega_{\mathrm{FS}}.
\]
On the standard affine chart centered at
\(x_0=([1:0],\ldots,[1:0])\), put
\[
 \rho_0(z)=\sum_{j=1}^n\log(1+|z_j|^2),\qquad
 \beta=\ddc|z|^2,\qquad Q=|z|^2-\rho_0.
\]
Since \(Q=O(|z|^4)\) with
\(|D^jQ(z)|\leq C|z|^{4-j}\) for \(j\leq2\), a cutoff at scale \(R\)
produces a smooth function \(\eta_R\), supported in \(B_{2R}\), such that
\[
 \|\ddc\eta_R\|_{\omega_0}\leq CR^2,\qquad
 \eta_R=Q\quad\text{on }B_R.
\]
For \(R\) small,
\[
 \omega:=\omega_0+\ddc\eta_R
\]
is K\"ahler, cohomologous to \(\omega_0\), and satisfies
\[
 \omega=\beta\quad\text{on }B_R.
\]

\medskip
\noindent\emph{Step 2: an \(m\)-subharmonic pole.}
Choose \(\chi\in C_c^\infty(B_R)\), with
\(0\leq\chi\leq1\) and \(\chi=1\) on \(B_{R/2}\), and set \(t=|z|^2\).  For
\(0<\varepsilon<R/4\), define
\[
 U_\varepsilon
 =-\delta\chi(z)(t+\varepsilon^2)^{-\lambda}.
\]
On the annulus \(B_R\setminus B_{R/2}\), the functions
\(\chi(t+\varepsilon^2)^{-\lambda}\) have uniformly bounded \(C^2\)-norm.
Thus one can choose \(\delta>0\), independently of \(\varepsilon\), so that
\(\omega+\ddc U_\varepsilon\) is K\"ahler there.

On \(B_{R/2}\), let
\[
 L_\varepsilon
 :=\delta\lambda(t+\varepsilon^2)^{-\lambda-1}.
\]
Relative to \(\beta=\omega\), the eigenvalues of \(\ddc U_\varepsilon\)
are \(L_\varepsilon\) in the \(n-1\) tangential directions and
\[
 L_\varepsilon\frac{\varepsilon^2-\lambda t}{t+\varepsilon^2}
 \geq-\lambda L_\varepsilon
\]
in the radial direction.  Hence, for \(1\leq k\leq m\),
\[
 \sigma_k(\omega+\ddc U_\varepsilon)
 \geq
 \binom{n-1}{k}(1+L_\varepsilon)^{k-1}
 \left[
  \frac{n}{n-k}
  +\frac{n-k-k\lambda}{n-k}L_\varepsilon
 \right]>0,
\]
because
\(\lambda<(n-m)/m\leq(n-k)/k\).
Therefore \(U_\varepsilon\in\SH_m(X,\omega)\), and
\[
 U_\varepsilon\searrow
 U:=-\delta\chi(z)|z|^{-2\lambda},\qquad \sup_XU=0.
\]
In particular, \(U\in\SH_m(X,\omega)\) and has exactly one pole, at \(x_0\).

\medskip
\noindent\emph{Step 3: the measure and the sharp thresholds.}
On \(B_{R/2}\setminus\{0\}\), put
\(L=\delta\lambda t^{-\lambda-1}\).  The eigenvalues of
\(\omega+\ddc U\) are
\[
 1+L,\ldots,1+L,\ 1-\lambda L.
\]
Since $\gamma^m\wedge\beta^{n-m}=\binom{n}{m}^{-1}
\sigma_m(\lambda_1,\ldots,\lambda_n)\,\beta^n$ for a real $(1,1)$-form
$\gamma$ with $\beta$-eigenvalues $\lambda_1,\ldots,\lambda_n$, this gives
the exact formula
\begin{equation}\label{eq:model-hessian-density}
 H_m(U)
 =
 (1+L)^{m-1}
 \left(1+\frac{n-m-m\lambda}{n}L\right)\beta^n.
\end{equation}
Consequently,
\begin{equation}\label{eq:model-density-two-sided}
 H_m(U)\asymp |z|^{-2m(\lambda+1)}\omega^n
 \quad\text{near }x_0.
\end{equation}
The density is locally integrable precisely when
\(m(\lambda+1)<n\), which is our upper bound on \(\lambda\).

Let $F$ denote the smooth density of the classical Hessian measure of
$U$ on $X\setminus\{x_0\}$, so that there $H_m(U)=F\omega^n$.  On
$B_{R/2}$, where $\omega=\beta$, the same computation for $U_\varepsilon$,
whose radial eigenvalue is
$L_\varepsilon(\varepsilon^2-\lambda t)/(t+\varepsilon^2)$, gives
\begin{equation}\label{eq:model-regularized-density}
 \frac{H_m(U_\varepsilon)}{\beta^n}
 =(1+L_\varepsilon)^{m-1}
 \left[
  1+L_\varepsilon\,
  \frac{n\varepsilon^2+(n-m-m\lambda)t}{n(t+\varepsilon^2)}
 \right].
\end{equation}
Since $n-m-m\lambda>0$, the coefficient of $L_\varepsilon$ inside the
brackets lies in $[0,1]$; moreover $L_\varepsilon\leq\delta\lambda
t^{-\lambda-1}$.  Hence
\[
 0\leq\frac{H_m(U_\varepsilon)}{\beta^n}
 \leq(1+L_\varepsilon)^m
 \leq C\bigl(1+t^{-m(\lambda+1)}\bigr).
\]
The dominating function is integrable because $m(\lambda+1)<n$.  On every
compact subset of $X\setminus\{x_0\}$, the convergence
$U_\varepsilon\to U$ is smooth.  Consequently the densities of
$H_m(U_\varepsilon)$ converge to $F$ in $L^1(X,\omega^n)$, and hence
\begin{equation}\label{eq:total-variation-model}
 H_m(U_\varepsilon)\longrightarrow
 \mu:=F\omega^n
 \quad\text{in total variation}.
\end{equation}
Each smooth measure on the left has mass $V_\omega$ by Stokes' theorem, so
$\mu(X)=V_\omega$.  Since $F\in L^1(\omega^n)$, one has
$\mu(\{x_0\})=0$.

It remains to identify this density limit with the canonical Hessian measure
of the unbounded function $U$.  For $j\geq1$ set
\[
 U_j:=\max\{U,-j\}.
\]
Then $U_j\in\SH_m(X,\omega)\cap L^\infty(X)$.  On the open set
$\{U>-j\}$, the equality $U_j=U$ holds and $U$ is locally bounded.
Locality of bounded Bedford--Taylor products gives
\begin{equation}\label{eq:canonical-truncation-identification}
 \one_{\{U>-j\}}H_m(U_j)
 =\one_{\{U>-j\}}F\omega^n
 =\one_{\{U>-j\}}\mu.
\end{equation}
As $j\to\infty$, the right-hand side increases to $\mu$, because
$\{U=-\infty\}=\{x_0\}$ and $\mu(\{x_0\})=0$.  Hence
\[
 \int_{\{U>-j\}}H_m(U_j)
 =\mu(\{U>-j\})\longrightarrow V_\omega.
\]
Every bounded $U_j$ has total Hessian mass $V_\omega$, so equivalently
\[
 \int_{\{U\leq-j\}}H_m(U_j)\longrightarrow0.
\]
By the full-mass criterion and the canonical truncation definition in
\cite[Definition~5.1 and Lemma~5.2]{LuNguyen}, after the harmless volume
rescaling explained in Section~\ref{sec:finite-energy-solvability}, we
conclude that
\begin{equation}\label{eq:model-full-mass}
 U\in\mathcal E(X,\omega,m),\qquad H_m(U)=\mu.
\end{equation}
Thus there is no residual Hessian mass at the pole: the canonical
truncations recover exactly the measure in
\eqref{eq:total-variation-model}.

The singularity also gives, for \(s>0\),
\[
 \int_{B_{R/2}}|U|^s\omega^n
 \asymp\int_0^{R/2}r^{2n-1-2\lambda s}\,dr,
\]
and hence
\begin{equation}\label{eq:model-Ls-threshold}
 U\in L^s(X,\omega^n)
 \quad\text{if and only if}\quad s<\frac n\lambda.
\end{equation}
Using \eqref{eq:model-density-two-sided},
\[
 \int_{B_{R/2}}(-U)H_m(U)
 \asymp
 \int_0^{R/2}
 r^{2n-1-2m(\lambda+1)-2\lambda}\,dr,
\]
so that
\begin{equation}\label{eq:model-weighted-integral}
 \int_X(-U)\,H_m(U)<+\infty
 \quad\text{if and only if}\quad
 \lambda<\frac{n-m}{m+1};
\end{equation}
at equality the divergence is logarithmic.

We relate this integral directly to the energy $\mathcal E_m$ of
Section~\ref{sec:finite-energy-solvability}.  Put $\nu_j:=H_m(U_j)$.  By
\eqref{eq:canonical-truncation-identification} and since both $\nu_j$ and
$\mu$ have mass $V_\omega$, we have $\nu_j=\mu$ on $\{U>-j\}$ and
$\nu_j(\{U\leq-j\})=\mu(\{U\leq-j\})$.

If $\int_X(-U)\,d\mu=+\infty$, note that for a bounded nonpositive
$w\in\SH_m(X,\omega)$ every term in the definition of $\mathcal E_m(w)$ is
nonpositive, so $\int_X(-w)\,H_m(w)\leq(m+1)\bigl(-\mathcal E_m(w)\bigr)$.
Since $U\leq U_j$, the monotone definition of the energy gives
$\mathcal E_m(U)\leq\mathcal E_m(U_j)$, and therefore
\[
 \int_{\{U>-j\}}(-U)\,d\mu
 \leq\int_X(-U_j)\,d\nu_j
 \leq(m+1)\bigl(-\mathcal E_m(U)\bigr).
\]
Letting $j\to\infty$ and using $\mu(\{x_0\})=0$ gives
$\mathcal E_m(U)=-\infty$.

If $\int_X(-U)\,d\mu<+\infty$, then \eqref{eq:energy-comparison} applied to
$U_j$, together with
\[
 j\,\mu(\{U\leq-j\})\leq\int_{\{U\leq-j\}}(-U)\,d\mu,
\]
gives
\[
 -\mathcal E_m(U_j)
 \leq\int_X(-U_j)\,d\nu_j
 \leq\int_{\{U>-j\}}(-U)\,d\mu+j\,\mu(\{U\leq-j\})
 \leq\int_X(-U)\,d\mu,
\]
and the continuity of $\mathcal E_m$ along decreasing sequences
\cite[Section~6]{LuNguyen} gives $\mathcal E_m(U)>-\infty$.  Together with
\eqref{eq:model-full-mass} and \eqref{eq:model-weighted-integral}, this
proves
\begin{equation}\label{eq:model-energy-threshold}
 U\in\mathcal E^1(X,\omega,m)
 \quad\text{if and only if}\quad
 \lambda<\frac{n-m}{m+1}.
\end{equation}

\medskip
\noindent\emph{Step 4: top-degree realization.}
From \eqref{eq:model-density-two-sided},
\[
 F\in L^p(X,\omega^n)
 \quad\text{for every}\quad
 1\leq p<\frac{n}{m(\lambda+1)}.
\]
The upper endpoint is larger than $1$, so choose such a $p>1$.  The measure
$\mu=F\omega^n$ is nonnegative and has mass $V_\omega$.
Ko{\l}odziej's compact K\"ahler solvability theorem for $L^p$ densities
\cite{KolodziejActa}, together with the H\"older regularity theorem
\cite[Theorem~2.1]{KolodziejHolder}, yields a normalized
$\Phi\in\PSH(X,\omega)\cap C^{0,\alpha}(X)$ for some $\alpha>0$ such that
\[
 (\omega+\ddc\Phi)^n=\mu.
\]
By \cite[Proposition~3.6]{DinhNguyenCMAPH}, the classical super-potential
of \(\mu\) is H\"older continuous.

The measure $\mu=F\omega^n$ does not charge $m$-polar sets, since
these sets have zero $\omega^n$-volume.
Finally, \eqref{eq:model-full-mass} and the full-mass uniqueness theorem
\cite[Theorem~4.1]{DinewLu} show that \(U\) is the unique normalized full-mass solution
of \(H_m(V)=\mu\).  If
\[
 \frac{n-m}{m+1}\leq\lambda<\frac{n-m}{m},
\]
then \eqref{eq:model-energy-threshold} gives
\(U\notin\mathcal E^1(X,\omega,m)\).  Hence the equation has no
finite-energy solution.  Since every bounded $\omega$-$m$-subharmonic
potential has full Hessian mass and finite energy, it follows in particular
that the equation has no bounded solution.
This completes the proof of Proposition~\ref{prop:projective-obstruction}.

We summarize the meaning of the obstruction.
\begin{remark}
The measure \(\mu\) has a H\"older classical super-potential, but its
lower-order equation \(H_m(V)=\mu\) has no finite-energy solution.  Thus
top-degree super-potential regularity cannot replace the same-order Hessian
hypothesis in Theorem~\ref{thm:weighted-solvability}.
\end{remark}

\subsection{Open questions}
\label{sec:open-questions}

The results above lead to several questions.  They are stated for compact
K\"ahler manifolds and $1\leq m<n$, where the Hessian and Monge--Amp\`ere
test spaces are different.

\begin{enumerate}[label=\textup{(Q\arabic*)},leftmargin=*,itemsep=0.75\baselineskip]
\item \emph{A converse for Hessian integration functionals.}
Let $\mu$ be a positive measure of mass $V_\omega$ such that
\[
 \left|\int_X(u-v)\,d\mu\right|
 \leq C\|u-v\|_{L^1(\omega^n)}^\beta,
 \qquad u,v\in\cP_m(X,\omega).
\]
Does there exist $\varphi\in\SH_m(X,\omega)\cap C^{0,\gamma}(X)$, for
some $\gamma>0$, satisfying $\Hm(\varphi)=\mu$?  For $m=n$ this is the
Dinh--Nguyen characterization; the corresponding question is open for
$m<n$.

\item \emph{Optimal functional exponents.}
If $\varphi\in C^{0,\alpha}(X)\cap\SH_m(X,\omega)$, what is the largest
exponent $\theta=\theta(m,\alpha,n)$ for which
\[
 \int_X |u-v|\,\Hm(\varphi)
 \leq C\|u-v\|_{L^1(\omega^n)}^\theta
\]
holds on $\cP_m(X,\omega)$?  The exponent $(\alpha/2)^m$ in
Theorem~\ref{thm:main-necessary} is the output of the present induction;
its optimality is unknown.

\item \emph{Regularity from the functional condition.}
Which hypotheses on a Hessian integration functional imply a quantitative
domination by Hessian capacity?  Such a domination is the natural bridge from
the functional condition to boundedness, continuity, or H\"older continuity
of solutions of the Hessian equation.

\item \emph{Hermitian backgrounds.}
Can the forward estimate of Theorem~\ref{thm:main-necessary} be extended to
a compact Hermitian manifold?  The bounded Hessian calculus in that setting
is available in \cite{KolodziejNguyenWeak,KolodziejNguyen}, but the present
argument uses closedness both in the mixed-current mass identity and in
Bedford--Taylor integration by parts.  A replacement retaining a uniform
H\"older modulus would be needed.
\end{enumerate}

\section*{Authors' affiliations}

\noindent Nguyen Quang Dieu, Department of Mathematics, School of
Mathematics and Computer Science, Hanoi National University of Education,
Hanoi, Viet Nam.\\
Email: ngquang.dieu@hnue.edu.vn

\medskip
\noindent Pham Le Khanh (corresponding author), Department of
Mathematics, School of Mathematics and Computer Science, Hanoi National
University of Education, Hanoi, Viet Nam.\\
Email: khanhlepham2002@gmail.com

\section*{Disclosure statement}

No potential conflict of interest was reported by the author(s).

\section*{Funding}

[Funding information to be completed by the authors.]

\end{document}